\documentclass[12pt,a4paper]{amsart}

\usepackage{bbm}

\usepackage{hyperref}

\def\mmat #1,#2,#3,#4,{\text{\small\arraycolsep=3pt $
\begin{pmatrix}#1&#2\\#3&#4\end{pmatrix}$}}

\usepackage{DLdef1}

\newcommand\vect{\mathfrak{vect}}
\newcommand\gl{\mathfrak{gl}}
\newcommand\po{\mathfrak{po}}
\newcommand\Specf{\mathop{\mathrm{Specf}}\nolimits}
\newcommand\Lee{\mathbb{L}}
\newcommand\Sym{\mathop{\mathrm{Sym}}\nolimits}
\renewcommand\Char{\mathop{\mathrm{Char}}\nolimits}

\begin{document}

\title{On realisations of the Steenrod algebras}

\author{Alexei Lebedev${}^{\mathrm {a}}$, Dimitry Leites${}^{\mathrm{b}*}$}

\address{${}^{\mathrm{a}}$Equa
Simulation AB, R{\aa}sun\-da\-v\"agen 100, Solna, Sweden; alexeylalexeyl@mail.ru\\
${}^{\mathrm{b}}$
Department of mathematics, Albanov\"agen 28, SE-114 19, Stockholm,
 Sweden; dimleites@gmail.com\\
 ${}^*$Corresponding author}

\begin{abstract} 
The Steenrod algebra can not be realised as an enveloping of
any Lie superalgebra. We list several problems that suggest a need to modify
the definition of the enveloping algebra, for example, to get rid of certain
strange deformations which we qualify as an artefact of the inadequate
definition of the enveloping algebra in positive characteristic. P. Deligne
appended our paper with his comments, hints and open problems.
\end{abstract}

\thanks{We are thankful to P. Deligne for his help
and encouragement,  to the International Max Planck Research School and MPIMiS-Leipzig,
for financial support and most creative environment during 2004--6.}


\makeatletter
\@namedef{subjclassname@2020}{\textup{2020} Mathematics Subject Classification}
\makeatother

\maketitle

\thispagestyle{empty}

\section{Introduction}

This is an updated version of our note published in not easily available \textit{J.~Prime Research in Mathematics}, \textbf{2} no.~1 (2006), 1--13 and preprinted in \texttt{MPIMiS preprint 131/2006}; this depositary was disabled since the end of 2024.
 
Hereafter, $p>0$ is the characteristic of the ground field $\Kee $.

\subsection{Motivations} In the mid-1970s, Bukhshtaber and Shokurov~\cite{BS} interpreted
the Landwe\-ber-Novikov algebra as the universal enveloping algebra of
the Lie algebra of the vector fields on the line with coordinate $t$, and with the
coefficients of $\frac d{dt}$
vanishing at the origin together with their first derivative.
Therefore, when (at about the same time) P.~Deligne told one of us (DL) that
Grothendieck told him what sounded (to DL) like a similar interpretation of
the Steenrod algebra; it did not alert the listener, although one should be very
careful when $p > 0$. From that time on till recently, DL remembered Deligne's
information in the following form
\begin{equation}
\label{(1)}
\text{\parbox{.8\textwidth}{%
\textsl{``The Steenrod algebra $\mathfrak A (2)$ is isomorphic to the enveloping algebra
$U(\mathfrak g)$ of a subsuperalgebra $\mathfrak g $ of the Lie superalgebra of contact vector
fields on the $1|1$-dimensional superline, whose generating functions
vanish at the origin together with their first derivative''}%
}}
\end{equation}
but the precise statement was never published and with time DL forgot what
at that time he thought he understood from Deligne what Grothendieck meant
under $\mathfrak g $. (For a precise statement in Deligne's own words, quite distinct from
\eqref{(1)}, see \S~4.) Somewhat later, Bukhshtaber~\cite{Bu} published a paper whose title
claims to interpret the Steenrod algebras $\mathfrak A (p)$ for $p > 2$ in terms similar to
\eqref{(1)}, namely as (isomorphic to)
\begin{equation}
\label{(2)}
\textsl{``the enveloping algebra of the supergroup of $p$-adic diffeomorphisms
of the line''.}%
\end{equation}
The body of the paper~\cite{Bu} clarifies its cryptic title (it is deciphered as meant
to be ``the enveloping algebra of a subsuperalgebra of the Lie superalgebra
of vector fields on the $1|1$-dimensional superline in characteristic $p > 2$''),
but nowhere actually states that the Steenrod algebra $\mathfrak A (p)$ is identified with
$U(\mathfrak g )$ for any $\mathfrak g $. Instead, $\mathfrak A (p)$ is realised by differential operators but no
description of the totality of these operators in more ``tangible'' terms, e.g.,
like the graphic (but wrong, as we will see) descriptions~\eqref{(1)}, \eqref{(2)} is offered; this
is an open (but perhaps unreasonable) problem, cf.~\cite{Wp} with \S~4.

Our initial intention was to explicitly describe the subsuperalgebra $\mathfrak g $ of the
Lie superalgebra of vector fields on the $1|1$-dimensional superline for which
$\mathfrak A (p) \simeq U(\mathfrak g )$ as we remembered~\eqref{(1)}; we also wanted to decipher~\eqref{(2)}. However,
having started, we have realised that we do not understand even what
$U(\mathfrak g )$ is if $p > 0$. More precisely, it is well-known (\cite{S}) that there are two
versions of the enveloping algebras (a ``usual'' one and a restricted one), but
it seems to us that there are many more versions. An \textbf{open problem} is to
give the appropriate definition of $U(\mathfrak g )$ (this is definitely possible, at least, for
classical Lie superalgebras with Cartan matrix) and related notions, such as
representations and (co)homology of $\mathfrak g $.

So we begin with a discussion of the notion of $U(\mathfrak g )$, and next pass to
realisations of the Steenrod algebras. We conclude that, under conventional
definitions (\cite{S}), there is no Lie superalgebra $\mathfrak g $ such that $\mathfrak A (p) \simeq U(\mathfrak g )$.

This result is not appealing: we hoped to clarify known realisations, not
make a negative statement that $\mathfrak A (p)$ is NOT something. There are, however,
realisations of $\mathfrak A (p)$ by differential operators (\cite{Bu,Wd}). These realisations,
although accepted, still look somewhat mysterious to us. In his comments (see \S~4), Deligne suggests a positive characterisation of~ $\mathfrak A (p)$.

\subsection{Notations} Let $T(V )$ be the tensor algebra of the superspace $V $, let
$S(V )$ and $\Lambda (V )$ be the symmetric and exterior algebra of the space $V$ , respectively.
For a set $x = (x_1,\ldots,x_n)$ of indeterminates that span $V$ , we write $T[x]$
or $S[x]$ or $\Lambda [x]$, respectively. Let $\mathbf Z _+$ be the set of nonnegative integers.

As an abstract algebra, the Steenrod algebra $\mathfrak A (p)$ is defined  as follows ($\deg \beta  = 1$):
\begin{equation}
\label{(3)}
\mathfrak A (p)=\begin{cases}
\{T[P^i] \mid \deg P^i=2i(p-1)\quad \text{for}\quad i\in \mathbf Z _+]\otimes \Lambda [\beta ])/I(p)\}&\text{for }p>2,\\
\{T[P^i] \mid \deg P^i=i \quad \text{for}\quad i\in \mathbf Z _+])/I(2)\}&\text{for }p=2,
\end{cases}
\end{equation}
where the ideal of relations $I(p)$ for $p > 2$ is generated by the \textit{Adem relations}
\newlength\flap
\newlength\flapm
\newlength\flapmo
\settowidth\flap{$\scriptstyle [a/p]$}
\settowidth\flapm{$\scriptstyle \rbrack (a-1)/p]$}
\flapmo=\flapm
\addtolength\flapm{-\flap}
\begin{equation}
\label{(4)}
\begin{alignedat}{2}
P^aP^b={}&\sum\limits _{i=1}^{[a/p]}(-1)^{a+i}\binom {(p-1)(b-i)-1}{a-pi}P^{a+b-i}P^i&\quad &\text{for }a<pb,\\
P^a \beta P^b={}&\sum\limits _{i=0}{[a/p]}(-1)^{a+i}\binom {(p-1)(b-i)}{a-pu}\beta P^{a+b-i}P^i-{}\\
&\sum\limits _{i=0}^{\hspace*{-.5\flapmo}\lefteqn{\scriptstyle [(a-1)/p]}\phantom {\scriptstyle [a/p]}}\hspace*{.5\flapm}(-1)^{a+i}\binom {(p-1)(b-i)-1}{a-pi-1}P^{a+b-i}\beta P^i&\quad &\text{for }a \leqslant pb,
\end{alignedat}
\end{equation}
whereas $I(2)$ is generated by
\begin{equation}
\label{(5)}
P^aP^b=\sum\limits _{i=1}^{[a/2]}(-1)^{a+i}\binom {b-i-1}{a-2i}P^{a+b-i}P^i \quad \text{for}\quad a<2b.
\end{equation}
\textbf{Remark.} For $p = 2$, the $P^i$ are usually denoted $Sq^i$.

\subsection{Lie superalgebras for $p = 2$} Observe that, for $p \ne 2$, for any Lie
superalgebra $\mathfrak g $ and any odd $x$, we have
\begin{equation}
\label{(6)}
[x, x] = 2x^2.
\end{equation}
In other words, there is a \textit{squaring operation}
\begin{equation}
\label{(7)}
x^2 =\frac 12\,
 [x, x]
\end{equation}
and to define the bracket of odd elements is the same as to define the squaring,
since
\begin{equation}
\label{(8)}
[x, y] = (x + y)^2 - x^2 - y^2 \quad \text{for any}\quad  x, y\in \mathfrak g _{\overline 1}.
\end{equation}

A \textit{Lie superalgebra} for $p = 2$ is a superspace $\mathfrak g $ such that $\mathfrak g _{\overline 0}$
is a Lie algebra, $\mathfrak g _{\overline 1}$
is an $\mathfrak g _{\overline 0}$-module (made into the two-sided one by symmetry) and on $\mathfrak g _{\overline 1}$
a
\textit{squaring} (roughly speaking, the halved bracket) is defined
\begin{equation}
\label{(9)}
\begin{aligned}
&x \mapsto x^2 \text{ such that }(ax)^2 = a^2x^2 \text{ for any }x\in \mathfrak g _{\overline 1} \text{ and }a\in \Kee;\\
&\text{and the map }(x, y)\mapsto (x + y)^2 - x^2 - y^2 \text{ is bilinear}\\
&\text{and }\mathfrak g _{\overline 0}\text{-invariant, i.e., } [x, y^2] = (\ad_y)^2(x)\text{ for any }x \in  \mathfrak g _{\overline 0}\text{
and }y \in  \mathfrak g _{\overline 1}.
\end{aligned}
\end{equation}
Then,  the bracket (i.e., product in $\mathfrak g $) of odd elements is defined to be
\begin{equation}
\label{(10)}
[x, y] := (x + y)^2 - x^2 - y^2.
\end{equation}
The Jacobi identity for three odd elements is replaced by the following relation:
\begin{equation}
\label{(11)}
[x, x^2] = 0 \quad \text{for any}\quad  x \in  \mathfrak g _{\overline 1}.
\end{equation}
This completes the definition unless the ground field is $\Zee /2$ in which case we
have to add the condition
\begin{equation}
[x, y^2] = (\ad_y)^2(x)\quad \text{for any}\quad  x \in \mathfrak g \quad \text{and}\quad  y \in  \mathfrak g _{\overline 1}
\label{(12)}
\end{equation}
which makes \eqref{(11)} and the last line in \eqref{(9)} redundant and replaces them over
any field. The \textit{restricted} Lie superalgebras are classically defined; for various versions of restrictedness in the case where $p=2$, see \cite{BLLS1}.

\subsection{Divided powers} For $p > 0$, there are many analogs of the polynomial
algebra. These analogs break into the two types: the infinite dimensional
ones and finite dimensional ones. The \textit{divided power algebra} in indeterminates
$x_1,\ldots , x_m$ is the algebra of polynomials in these indeterminates, so, as space,
it is
$$
\mathcal O (m) = \Span \{x_1^{(r_1)}\ldots x_m^{(r_m)}\mid r_1,\ldots , rm \geqslant  0\}
$$
with the following multiplication:
$$
(x_1^{(r_1)}\ldots x_m^{(r_m)})\cdot (x_1^{(s_1)}\ldots x_m^{(s_m)}):=\prod\limits _{i=1}^m \left (\dbinom {r_i+s_i}{r_i}\right ) x_i^{(r_i+s_i)}.
$$

For a \textit{shearing vector of heights of the indeterminates} $\underline  N = (N_1,\ldots,N_m)$, set
$$
\mathcal O (m,\underline  N)\quad  (\text{or } \Kee [u;\underline  N]) := \Span\{x_1^{(r_1)}\ldots x_m^{(r_m)}\mid 0 \leqslant r_i<p^{N_i},i=1,\ldots ,m\},
$$
where $p^\infty  = \infty $. If $N_i < \infty $ for all $i$, then $\dim \mathcal O (m,\underline N) < \infty $.

Observe that only the conventional polynomial algebra and the divided power one with
$\underline N = (1,\ldots , 1)$ are generated by the indeterminates that enter its definition.
For any other value of $\underline N$, we have to add $x_i^{(p^{k_i })}$ for every $k_i$ such that $1 <
k_i < N_i$.

If an indeterminate $x$ is odd, then the corresponding height of shearing vector is
equal to 1 for $p = 2$; for $p > 2$, we postulate
$$
x^2 = 0 \text{ and anticommutativity of odd elements.}
$$

\section{The enveloping algebras of Lie algebras for $p > 0$}

It looks strange that the following problem was never discussed in the literature.
For $p > 0$, it seems natural --- in view of the Poincar\'e--Birkhof--Witt
theorem --- to have as many types of universal enveloping algebras, as there
are analogs of symmetric algebras or algebras of divided powers.

Of the variety of such hypothetical definitions of enveloping algebras (the
usual one and the ones labelled by various values of the shearing vector
$\underline N$), only two are being considered: the usual $U(\mathfrak g )$ and the one corresponding
to $\underline N = (1,\ldots , 1)$.

We hope that there exist $U(\mathfrak g ;\underline N)$ --- analogs of $U(\mathfrak g )$, such that $\gr U(\mathfrak g ;\underline N) \simeq \mathcal O (\dim \mathfrak g ;\underline N)$. Such $\underline  N$-dependent definitions of $U(\mathfrak g )$ do exist, at least for simple
finite dimensional complex Lie algebras, see~\cite{St} (take the $\Zee $-form of $U(\mathfrak g )$
described in~\cite{St} and tensor by $\Kee $) and we hope that the \textbf{open problem} to
generalize the definition to arbitrary algebras $\mathfrak g $ is not difficult to solve. We
suspect that these $U(\mathfrak g ;\underline N)$ were ignored because they are not generated by the
initially declared indeterminates (or the space $\mathfrak g $ they span) and the necessity
to add extra generators (depending on $\underline N$) was too unusual: to preserve a
one-to-one-correspondence between representations of $U(\mathfrak g )$ and representations  of $\mathfrak g $, we have to amend the definitions.

The notion of the universal enveloping algebra is motivated, first of all,
by the representation theory. So let us give more reasons, other than the
PBW theorem, to consider the non-conventional universal enveloping algebras
corresponding to any value of the shearing vector~$\underline N$.

\subsection{The induced and coinduced modules} Since any derivation of a
given algebra is completely determined by its values on every generator of
the algebra, the Lie algebra of all derivations of $\Kee [u;\underline N]$ is much larger than
the Lie algebra of \textit{distinguished}\footnote{The term \textit{special derivation}, often used, causes confusion with divergence-free derivations.} derivations, where one partial derivative replaces several conventional partial derivatives over divided powers of $u_j$:
\begin{equation}
\partial _i(u_j^{(k)})=\delta _{ij}u_j^{(k-1)}.
\label{(13)}
\end{equation}
In what follows, speaking about Lie algebras of vector fields (briefly: vectorial
algebras) we consider only distinguished derivatives, e.g., in \eqref{(14)}.

The simple vectorial Lie algebra for $p = 0$ has only one parameter: the
number of indeterminates. If $\Char \Kee  = p > 0$, the vectorial Lie algebras
acquire one more parameter: $\underline N$. For Lie superalgebras, $\underline N$ concerns only the
even indeterminates. Let
\begin{equation}
\mathfrak {vect}(m;\underline  N |n) \text{ a.k.a } W(m;\underline  N|n) :=\mathfrak { der}\Kee [u;\underline N]
\label{(14)}
\end{equation}
be the general vectorial Lie algebra.

The induced and coinduced modules are natural classes of modules over Lie
algebras. Over $\Cee $, the modules of (formal) tensor fields constitute a natural
class of modules. In particular, a most natural --- $(x)$-adic --- filtration in the
polynomial algebra $\Kee [x_1,\ldots , x_m]$, induces a filtration (\textit{Weisfeiler filtration})
in the Lie algebra $\overline {\mathfrak { vect}(m)} :=\mathfrak { der}(\Kee [[x]])$. The associated grading is given by
setting $\deg x_i = 1$ for all $i$. Let
$$
\mathcal L  =\mathcal L_{-1}\subset \mathcal L_0 \subset \mathcal L_1 \subset \ldots 
$$
be the Weisfeiler filtration of $\mathfrak L  :=\mathfrak {vect}(m)$; let $L_i :=\mathcal  L_i/\mathcal L_{i+1}$.

Let $V$ be a $\mathfrak {gl}(m)$-module, considered as a $\mathcal L_0$-module such that $\mathcal L_iV = 0$
for $i > 0$. We define the induced and coinduced modules over $\mathcal L $ as
\begin{equation}
\Ind\nolimits ^{\mathcal L}_{\mathcal L_0}(V ) = U(\mathcal L )\otimes _{U(\mathcal L_0)} V,\qquad  \Coind \mathcal  L_{\mathcal  L_0}(V ) = \Hom\nolimits_{U(\mathcal L_0)}(U(\mathcal L ), V ).
\label{(15)}
\end{equation}
In particular, the spaces of tensor fields of type $V$ (with fiber $V $) are coinduced $\mathcal 
L$-modules.

The spaces $\mathcal O (m;\underline N)$ and $\vect(m;\underline N)$ are $\vect(m;\underline  M)$-modules (well-defined
if $M_i \leqslant N_i$ for every $i$), but if we want to consider them as coinduced modules,
we need the unconventional universal enveloping algebras, namely for the
commutative Lie (super)algebra $\mathcal L _{-1}/\mathcal L _0 \simeq \Span(\partial _1,\ldots \partial _m)$, we define
\begin{equation}
U(\mathcal L _{-1}/\mathcal L_0;\underline N) := \mathcal O (m;\underline N). 
\label{(16)}
\end{equation}

The computation of deformations of $\vect(m;\underline N)$ is currently performed either
by painstaking calculations (\cite{DK,Dz}) or with the help of computer, see \cite{BGL}, but
both with the same --- conventional --- definition of $U(\mathfrak g )$. The late Grozman 
used his remarkable (highly appreciated by Etingof and his MIT students) \textit{SuperLie} package, now maintained by Krutov, see \cite{Gr},  to verify the rigidity of $\vect(m;\underline N)$ for small $m$ and
$\underline N = (1,\ldots , 1)$ for $p = 3$, see \cite{GL}, whereas, for $\underline N \ne (1,\ldots , 1)$, Dzhumadildaev and
Kostrikin \cite{DK} found lots and lots of infinitesimal deformations (2-cocycles)
all of which are mysterious; see also \cite{BLLS} on semi-trivial deformations in any characteristic.

On the other hand, recall that, for $p = 0$ and any $\mathfrak h $-module $M$, we have
(\cite{Fu})
\begin{equation}
H^q(\mathfrak g ; \Coind^{\mathfrak g }_{
\mathfrak h }(M)) \simeq H^q(\mathfrak h ;M);\quad H_q(\mathfrak g ; \Ind^{\mathfrak g }_{
\mathfrak h }(M)) \simeq  H_q(\mathfrak h ;M).
\label{(17)}
\end{equation}
Now, observe that, over vectorial Lie superalgebras $\mathfrak g $, the modules of tensor
fields are precisely the coinduced ones: $T(M) := \Coind^{\mathfrak g }_{\mathfrak g \geqslant 0}(M)$, where $\mathfrak g _{\geqslant  0} :=\underset {i \geqslant 0}{\oplus }\mathfrak g _i$
and $M$ is any $\mathfrak g _{\geqslant 0}$-module such that $\mathfrak g _{>0}M = 0$ for $\mathfrak g _{>0} =\underset {i>0}\oplus 
\mathfrak g _ i$, i.e., $M$ is,
actually, a $\mathfrak g _ 0$-module. In particular, let $\id_{\gl(m)}$ be the tautological (identity) $\gl(m)$-module.
Then, 
$$
W(m;\underline N) \simeq \mathcal O (m;\underline N) \otimes \id_{\gl(m)}
$$
which is a coinduced module if we define $U(\mathcal L _{-1}/\mathcal L _0;\underline N)$ my means of \eqref{(16)}.

So, for a conjectural cohomology theory $\text{``$H$''}^._{\underline N}$, the following rigidity conjecture
would be an corollary of the general theorem \eqref{(17)} and the mysterious infinitesimal deformations found in \cite{DK,Dz} should be considered as ``artefacts''
(except, perhaps, certain values of $p$ and $m$ for which $\text{``$H$''}^2(\gl(m); \id_{\gl(m)}) \ne  0$)
because $\text{``$H$''}^i_{
\underline B}(\gl(m); \id_{\gl(m)}) = H^i(\gl(m); \id_{\gl(m)})$ for $i$ small.

\subsubsection{Corollary \normalfont{(Conjecture)}} \textit{We have}
$$
\text{``$H$''}^2_{
\underline N}(W(m;\underline N);W(m;\underline N))\simeq \text{``$H$''}^2_{
\underline N}(\gl(m); \id_{\gl(m)}) = 0.
$$

\subsection{How to quantize?} The Poisson Lie (super)algebra $\po(2n|m)$ realised
on polynomials admits only one deformation as a Lie (super)algebra for $p=0$, see \cite{LS}.
After Dirac, physicists interpret this deformation as quantisation. Quantisation
deforms $\po(0|2m)$ into $\gl(\Lambda (m))$ and $\po(0|2m-1)$ into $\fq(\Lambda (m))$, see \cite{Lo}. What is the analog of this statement for
$p > 0$ and $\po(2n;\underline N|2m)$? The answer depends on how we understand $U(\mathfrak g )$.

\section{Main result}

\ssbegin{Theorem}
For any $\Zee $-graded Lie superalgebra $\mathfrak g  ={\oplus  }\mathfrak g _ k$, we consider  the induced $\Zee$-grading $U(\mathfrak g ) ={\oplus }U(\mathfrak g )_k$. 

For any $p > 0$, if $P^0$ (or $Sq^0$) is a
scalar, there is no grading preserving isomorphism $f : \mathfrak A (p)\tto U(\mathfrak g )$ between
the Steenrod algebra $\mathfrak A (p)$ and the (common or restricted) universal enveloping
algebra of any $\Zee$-graded Lie superalgebra $\mathfrak g $ with the parity of elements of $\mathfrak g _ k$
being the same as that of $k$.\end{Theorem}

\begin{proof}
Suppose that such an isomorphism exists. First, let us show that
$\dim \mathfrak g _k$ is uniquely determined by the information on $\dim \mathfrak A (p)_i$. Let
$$
\mathfrak G _ k :=\bigoplus\limits _{i \leqslant k}
\mathfrak g _k
$$
as a ($\Zee$-graded) linear superspace. Clearly, $\dim \mathfrak G _ 0 = 0$ (since $\dim \mathfrak A (p)_i = 0$ for
$i < 0$). So, according to PBW theorem, $\dim U(\mathfrak g )_k = \dim \mathfrak g _k + d_k$, where $d_k$ is
equal to the dimension of the space of (super)symmetric polynomials on $\mathfrak G _{k-1}$
of weight $k$, if we consider the non-restricted (common) universal enveloping
algebra, or to the dimension of the space of (super)symmetric polynomials on $\mathfrak G _{k-1}$
of weight $k$ and degree $< p$ w.r.t. any even basic element, if we consider
the restricted universal enveloping algebra.

Since $\dim U(\mathfrak g )_k = \dim \mathfrak A (p)_k$, and $d_k$ is determined by dimensions and
parities of $\mathfrak g _i$ for $i < k$, one can find $\dim \mathfrak g _k$ for any $k$ by induction. The
following table illustrates this for $p = 2$, the non-restricted algebra $U(\mathfrak g )$ and small
values of $k$. In the table, the first row contains $k$; the second row contains
bases of $\mathfrak A (2)_k$; the third row contains bases of the spaces of (super)symmetric
polynomials on $\mathfrak G _{k-1}$ of degree $k$, where $L_i$ denotes a non-zero element of $\mathfrak g _k$;
the fourth row contains $\dim \mathfrak g _k$.
\begin{center}
\begin{tabular}{|c|c|c|c|c|c|c|c|}
\hline
1& 2& 3 &4 &5 &6 &7 &8\\
\hline
$Sq^1$&$ Sq^2$&$\begin{array}{@{}c@{}}
 Sq^3\\
Sq^2Sq^1
\end{array}
$&$\begin{array}{@{}c@{}}
Sq^4\\
Sq^3Sq^1
\end{array}
$&$\begin{array}{@{}c@{}}
Sq^5\\
Sq^4Sq^1
\end{array}
$&$\begin{array}{@{}c@{}}
Sq^6\\
Sq^5Sq^1\\
Sq^4Sq^2
\end{array}
$&$\begin{array}{@{}c@{}}
Sq^7\\
Sq^6Sq^1\\
Sq^5Sq^2\\
Sq^4Sq^2Sq^1
\end{array}
$&$\begin{array}{@{}c@{}}
Sq^8\\
Sq^7Sq^1\\
Sq^6Sq^2\\
Sq^5Sq^2Sq^1
\end{array}
$\\
\hline
---&---&
$L_2L_1$&$\begin{array}{@{}c@{}}
L^2_2\\
L_3L_1
\end{array}
$&$\begin{array}{@{}c@{}}
L_2^2
L_1\\
L_2L_3
\end{array}
$&$\begin{array}{@{}c@{}}
L^3_2\\
L_2L_3L_1
\end{array}
$&$\begin{array}{@{}c@{}}
L_6L_1\\
L_2^2
L_3\\
L^3_2
L_1
\end{array}
$&$\begin{array}{@{}c@{}}
L_6L_2\\
L^4_2\\
L_2^2
L_3L_1\\
L_7L_1
\end{array}
$\\
\hline
1 &1 &1 &0 &0 &1 &1 &0\\
\hline
\end{tabular}
\end{center}

Now, let us consider the case $p = 2$. From the table and similar computations
for a hypothetical restricted algebra $\mathfrak g $, we see that $L_1 = Sq^1$, $L_2 = Sq^2$
are elements of $\mathfrak g $. If $Sq^0 = 0$, then $L_1L_2 = 0$, which can not be true. If $Sq^0$ is
a non-zero scalar, then, up to a non-zero scalar factor,
$$
(L_2)^2 = Sq^3Sq^1 \ne 0;\qquad  (L_2)^2L_1 = 0.
$$
This can hold only if we consider a restricted algebra, and $Sq^3Sq^1$ is an element
of $\mathfrak g $, proportional to $Sq^1$ --- which can not be true, since these two non-zero
elements have different weights.

Now, we consider the case $p > 2$. The computation of dimensions similar
to the above ones shows that the minimal weights in which $\mathfrak g $ has non-zero
elements are:
\[
\begin{cases}
\text{$1$, $2(p - 1)$, $2p - 1$, $2p^2 - 2$, $2p^2 - 1$, $2p^3 - 2$}&\text{for $\mathfrak g $ non-restricted},\\
\text{$1$, $2(p - 1)$, $2p - 1$, $2(p - 1)p$, $2p^2 - 2$, $2p^2 - 1$, $2(p - 1)p^2$}&\text{for $\mathfrak g$  restricted}.
\end{cases}
\]

Since $\dim \mathfrak A (p)_{2(p-1)} = 1$, we see that $L_{2(p-1)} = P^1$ is an element of $\mathfrak g $. If
$P^0 = 0$, then, $(L_{2(p-1)})^2 = (P^1)^2 = 0$, which is false. It follows from the Adem
relations that
$$
(L_{2(p-1)})^p = (P^1)^p = 0,
$$
which can hold only in a restricted algebra. Then,  since $\dim \mathfrak A (p)_{2(p-1)p} = 1$,
we see that $L_{2(p-1)p} = P^p$ is an element of $\mathfrak g $. If $P^0 = \lambda \ne 0$, it follows from
the Adem relations that
\begin{align*}
&[L_{2(p-1)},L_{2(p-1)p}] = [P^1,P^p] = P^1P^p - P^pP^1 =\lambda P^{p+1} - P^pP^1;\\
&[L_{2(p-1)}, [L_{2(p-1)},L_{2(p-1)p}]] = [P^1, [P^1, P^p]] {}\\
&=P^1(\lambda P^{p+1} - P^pP^1) - (\lambda P^{p+1} - P^pP^1)P^1 {}\\
&=2(\lambda ^2P^{p+2} -\lambda P^{p+1}P^1 +\lambda P^pP^2) = 2 \lambda P^{p+1}P^1 \ne 0.
\end{align*}

The last expression must be an element of $\mathfrak g $, but $\mathfrak g $ does not have non-zero
elements of weight $2(p - 1)(p + 2)$, so we get a contradiction.
\end{proof}

\section{Pierre Deligne's comments in a letter to DL, May 23, 2006}

About Steenrod. What Grothendieck saw is the following (for $p$ odd).

\medskip

{\large \texttt{a)}} Quillen \cite{Q1}: for complex cobordism, one has
$$
\Omega U(B\Cee^\times)=\Omega U(Pt)[\eta ],\quad \text{where}\quad \deg(\eta ) = 2
$$
(with $B\Cee = \Pee^\infty (\Cee)$), and the group law of $B\Cee^\times$, deduced from that of $\Cee^\times$,
induces on $\Spec \Omega U(B\Cee^\times)$ a structure of formal group over $\Spec \Omega U(Pt)$. This
turns $\Spec \Omega U(Pt)$ into the scheme of formal group laws on the pointed formal
disc $\Specf(\Zee[\![t]\!])$:
$$
\Omega U(Pt) = \Zee[a_{i,j}\mid  i, j \geqslant 0, i + j > 0]/\text{identities},
$$
the identities expressing that $F(t, u) =\sum\limits a_{ij}t^iu^j$ is a formal group law. The
group scheme of automorphisms of the pointed formal disc hence acts by transport
of structures on $\Omega U(Pt)$. It is the group of
\begin{equation}
t \longmapsto \sum\limits a_it^i \qquad \text{($i \geqslant 1$, $a_1$ invertible)}.
\label{(18)}
\end{equation}
The action of the subgroup $\Gee_m\colon  t \mapsto at$ gives the half degree. If we consider
the subgroup with $a_1 = 1$, this action extends to a functorial action
on $\Omega U(X)$, compatible with products (Landweber operations, see~\cite{BS}). The
group scheme of transformations \eqref{(18)} has a double covering, with coordinates $\sqrt {
a_1}$ and the $a_i$ ($i \geqslant 2$). This double covering is again a group scheme, and it
contains the $\Gee_m$-subgroup ``$a_i = 0$ for $i \geqslant 2$'' (coordinate $\sqrt {a_1}$). The action of
this $\Gee_m$ gives the degree.

\medskip

{\large\texttt{b)}} This suggests that for any commutative ring $R$, and any 1-dimensional
formal group $G$ over $R$, possibly given with a trivialisation of its Lie algebra:
$\Lie(G) \overset \sim \tto R$, there could be a corresponding cohomology theory, functorial
in $G$. If $t$ is a parameter for $G$ (compatible with the trivialisation of the Lie
algebra), $G$ is given by
$$
\Omega 
U(Pt) \to  R
$$
and the theory would be obtained from complex cobordism by some ``derived
extension of scalars'', while the Landweber operations would ensure that the
result is independent of the choice of $t$, up to unique isomorphisms.

I am rather naive here; we are playing with (ringed) spectrum, not with
rings and their derived categories. I don't know what has been done, but
results are known: As I remember being told, the case where $\Spec(R)$ is a
complete intersection in $\Spec(\Omega U(P^t))$ is OK. This allows for the construction
of Morava's $K$-theories using this philosophy.

\medskip

{\large\texttt{c1)}} $G = \Gee_m^\wedge/\Zee$:
in each characteristic, we are in the open orbit of the action
on $\Spec \Omega U(Pt)$, so that the extension of scalars to $\Zee$ is an exact functor, and
one gets $K$-theory (Conner and Floyd \cite{CF}).

\medskip

{\large\texttt{c2)}} For a formal group over a field of char $p$, the (geometric) invariant is the
height, and one gets the Morava $K$-theories (\cite{DMea}).

\medskip

{\large\texttt{c3)}} For $\Gee_a^\wedge /\Fee_p$, one gets the ordinary mod $p$ cohomology. The group scheme
of automorphisms of $\Gee^\wedge _a$
(which are 1 on the Lie algebra) should hence act. It
is the group scheme
$$
A =\left \{t \longmapsto \sum\limits b_it^{p^i},
b_0 = 1 \right \},
$$
whose affine algebra is $\Zee[b_1,b_2,\ldots]$ (denoted $\mathcal O (A)$) the coproduct (giving the
group law) being defined by
\begin{gather*}
b_0 := 1,\\
\Delta b_k =\sum\limits _{\ell+m=k}b_\ell \otimes b_{k-\ell }^{p^\ell }.
\end{gather*}
As shown by Milnor, this group indeed acts functorially on $H^{\bcdot}(\cdot ; \Fee_p)$, see \cite
[Th. 3, page 162]{Mi}. ``Action'' means ``comodule structure $H^{\bcdot} \to H^{\bcdot}\otimes \mathcal O (A)$''.

This does not capture the odd part of the story, for which I lack understanding.
What Milnor says is that (for $p$ odd)
$$
\Spec(H^{\bcdot}(B \Zee/p,\Zee.p)),
$$
with the group law coming from that of $B\Zee/p$, is  $\Gee^+_a \times 
\Gee^-_a$ (where obviously, one factor is an ``even'' group, the other one is an ``odd'' one, representing,
respectively, the functors $C \mapsto C_\ev$
and $C \mapsto C_\od$
for any supercommutative super ring $C$), i.e., $H^.
=
\Zee/p [t,\tau]$ with $t$ even and $\tau $ odd, and the group law
$$
(t',\tau ')+(t'',\tau '')=(t'+t'',\tau '+\tau '').
$$
If $B$ is the super group scheme\footnote{For the definition imitating Grothendieck's definition of groups schemes (\cite{MAG}), see \cite{Lsos}.} of automorphisms
of $G = \Gee^+_a \times 
\Gee^-_a
$, respecting the filtration $\Lie\Gee^+_a
\subset \Lie\Gee$ and acting trivially
on the successive quotients, the group $B$ acts functorially on $H^{\bcdot}(X,\Zee/p)$,
respecting the cup-product [and one could add to it a $\Gee_m$ giving the degree].
The action on $H^{\bcdot}(B\Zee/p,\Zee/p)$ is the one defining $B$, and the affine algebra
$\mathcal O (B)$ is the dual of the Steenrod algebra \cite
[Th. 2, page 159]{Mi}.

I would hope that the odd part of the story is analogous to the following
fact: if $k$ is a quotient of a ring $R$, then $\Ext^i_
R(k, k)$ acts on $H^{\bcdot}(M {\overset \Lee \otimes }_Rk)$ for any
$M$ in $D^-(R)$.

\medskip

Other comments on the text with Lebedev.

\medskip

Other convenient definitions of the space of quadratic forms on a projective
module $M$:

\medskip

\textbullet{}
$\Sym^2M^\vee $, where $Sym^2 = \text{covariants}$ of $S_2$ acting on 2-tensors.\\
equivalently: the space of quadratic form is the cokernel of the map 
\[
C \mapsto C(X, Y ) - C(Y,X)
\]
on the space of bilinear maps,

\medskip

\textbullet{}
the dual of $\Gamma ^2(M)$ (divided power = symmetric 2-tensors)\\
--- If $G$ is a smooth algebraic group on $\Spec(R)$, a reasonable analog of what
$U(\mathfrak g )$ is in characteristic 0 is the algebra of left-invariant differential operators
on $G$. As a coalgebra, it is the dual of the completion of $G$ at the unit element.
It is, I think, what Dieudonn\'e calls the \textit{hyperalgebra of the group}. It cannot
be constructed from the Lie algebra. For instance, $\Lie\Gee_a = \Lie\Gee_m = R$, but
for $G_a$ one gets the $\frac {\partial ^i_x}{
i!}$, and for $\Gee_m$ the binomial $\binom {x \partial _x}i$, where the choice of
generator is crucial.

\subsection{Pierre Deligne's comments in a letter to DL, September 1, 2006}
As I am a geometer, groups are more congenial to me than Lie algebras, and
it does not bother me that in characteristic $p>0$ the Lie algebra of a group does
a poor job of controlling it. If I want to have all relevant ``divided powers''
for a given group, I just take as starting point the bialgebra of left invariant
differential operators. This ``is the same'' as giving the formal group ($\mathcal O (G^\wedge ) =\text{dual}$) and is, if I remember right, what Dieudonn\'e calls a hyper (Lie?) algebra.

Lie algebras with a $p$th power operation (= restricted), on the other hand,
are exactly the same things as algebraic groups equal to the Kernel of Frobenius.

So, I am more happy with Steenrod ``being'' a (super) group scheme than
it being some kind of enveloping algebra.

Even in the even case, characteristic 2 and 3 are tricky, and I am not sure
one definition is suitable for all applications.\footnote{Indeed! See \cite{BLLS1} on restrictedness and \cite{KLLS} on ``even rules", see \cite{Del}.}


\begin{thebibliography}{999990}
\bibitem[Bu]{Bu} 
Bukhshtaber V. The Steenrod algebra is the enveloping algebra of the supergroup
of $p$-adic direomorphisms of the line. In: \cite{SE}, 224--228pp.
\end{thebibliography}
\end{document}